\documentclass[12pt]{article}

\usepackage{amsfonts,amsmath,amssymb,latexsym,xcolor,mathrsfs,tikz,breqn,extarrows}
\usepackage{amsthm,authblk}
\usepackage{multirow}
\usepackage{amssymb}
\usepackage{amsmath}
\usepackage{amsthm}
\usepackage{tikz}
\usepackage{mathtools}
\usepackage{calc}
\usepackage{tikz}
\usepackage{pgfplots}
\usetikzlibrary{patterns,arrows,decorations.pathreplacing}
\usepackage{diagbox}
\usepackage{makecell}

\def\q{\hfill\rule{1ex}{1ex}}
\def\0{\emptyset}

\def\q{\hfill\rule{1ex}{1ex}}

\usepackage{amsthm}

\newtheorem{theorem}{Theorem}[section]
\newtheorem{definition}[theorem]{Definition}
\newtheorem{lemma}[theorem]{Lemma}

\newtheorem{prop}[theorem]{Proposition}

\newtheorem{conj}[theorem]{Conjecture}
\newtheorem{prob}[theorem]{Problem}

\usepackage{algorithm}
\usepackage{algpseudocode}
\usepackage{graphicx}
\usepackage{amsmath}
\usepackage{pstricks,pgf,subcaption,amssymb}
\begin{document}
\title{\bf The saturation number of wheels 
	 }
\author[]{
Yanzhe Qiu}
\author[2]{
Zhen He}
\author[]{
Mei Lu}
\author[]{
Yiduo Xu\thanks{Corresponding author. E-mail:\texttt{xyd23@mails.tsinghua.edu.cn}}\;}

\affil[]{\small Department of Mathematical Sciences, Tsinghua University, Beijing 100084, China}
\affil[2]{\small School of Mathematics and Statistics, Beijing Jiaotong University, Beijing 100044, China.}
\date{}

\maketitle\baselineskip 16.3pt

\begin{abstract}
	A graph $G$ is said to be $F$-free, if $G$ does not contain any copy of $F$. $G$ is said to be $F$-semi-saturated, if the addition of any nonedge $e \not \in E(G)$ would create a new copy of $F$ in $G+e$. $G$ is said to be $F$-saturated, if $G$ is $F$-free and $F$-semi-saturated. The saturation number $sat(n,F)$ (resp. semi-saturation number $ssat(n,F)$) is the minimum number of edges in an  $F$-saturated (resp. $F$-semi-saturated) graph of order $n$. In this paper we proved several results on the (semi)-saturation number of the wheel graph $W_k=K_1 \vee C_k$. Let $k,n$ be positive integers with $k \geq 8$ and $n \geq 56k^3$, we showed that $(s)sat(n,W_k)=n-1+(s)sat(n-1,C_k)$. We also establish the lower bound of semi-saturation number of $W_k$ with restriction on maximum degree.

	\end {abstract}
	
	{\bf Keywords.} saturation number, semi-saturation number, cycle, wheel, extremal graph
	
\section{Introduction}

In this paper we only consider finite, simple and undirected graphs. For a graph $G$, we use $V(G)$ to denote the vertex set of $G$, $E(G)$ the edge set of $G$, $|G|$ the order of $G$ and $e(G)$ the size of $G$. For distinct $V_1,V_2 \subseteq V(G)$, let $E(V_1,V_2)$ be the set of edges between $V_1,V_2$ and $e(V_1,V_2)=|E(V_1,V_2)|$. If $V_1=\{u\}$, then we abbreviate as $E(u,V_2)$ and $e(u,V_2)$. For a positive integer $k$, let $[k]:=\{1,2,\ldots,k\}$. Denote by $C_k$ the $k$-vertices cycle. For a graph $G$ and $u \in V(G)$, let $N_G(u)=\{ v : uv \in E(G) \}$ be the \textit{neighbourhood} of $u$, and $d_G(u)=|N_G(u)|$  the degree of $u$. We use $\delta(G)$ and $\Delta(G)$ to denote the minimum   and the maximum  degree of  $G$ respectively. Without confusion, we abbreviate as $N(u)$, $d(u)$, $\delta$ and $\Delta$, respectively. For a vertex set $A $, let $N(A)= \bigcup \limits_{u \in A} N(u) \, \setminus \, A$. Let $N[u]=N(u) \cup \{u\}$ and $N[A] = N(A) \cup A$.
 For graphs $G_1,G_2,\ldots,G_t$, let $G_1 \cup G_2 \cup \ldots \cup G_t$ be the union of vertex-disjoint copy of $G_1,G_2,\ldots,G_t$. If $G_1=G_2=\ldots=G_t$ then we abbreviate the union as $tG_1$. Let $G_1 \vee G_2$ be the join of $G_1$ and $G_2$, obtained by adding all edges between $G_1$ and $G_2$ in $G_1 \cup G_2$.  For a given graph $G$, we say $u \in V(G)$ is a \textit{conical vertex} if $d_G(u)=|G|-1$.  The \textit{distance} $d(u,v)$ between two vertices $u,v \in V(G)$ is the number of edges contained in the shortest path connecting $u$ and $v$.

	Given graphs $G$ and $F$, we say  $G$ is \textit{$F$-free}, if $G$ does not contain any copy of $F$. We say  $G$ is \textit{$F$-semi-saturated}, if the addition of any nonedge $e \not \in E(G)$ would create a new copy of $F$ in $G+e$. We say  $G$ is \textit{$F$-saturated}, if $G$ is $F$-free and $F$-semi-saturated. The \textit{saturation number} and \textit{semi-saturation number} of $F$ is denoted by
\begin{equation*}
\begin{aligned}
	   sat(n,F) & =\min \{ e(G):G  \; \text{is}  \; F\text{-saturated  and } |G|=n  \} \, , \\
	   ssat(n,F) & =\min \{ e(G):G  \; \text{is}  \; F\text{-semi-saturated  and } |G|=n  \} \, .
\end{aligned}
\end{equation*}
By the definition of (semi)-saturation number, it is clear that $sat(n,F) \geq ssat(n,F)$. The first (semi)-saturation problem was studied by Erd\H os, Hajnal and Moon \cite{EHM} in 1964, determining that $ssat(n,K_r)=sat(n,K_r)=(r-2)(n-1)-\frac{(r-2)(r-1)}{2}$. For readers interested in saturation problem, we refer to the survey \cite{survey}.

	 To compute the (semi)-saturation number of a fixed graph $F$ on $n$ vertices, an efficient way is to partition all $F$-(semi)-saturated graphs  into several groups based on some parameter. Let $\Delta=\Delta(G)$ be the maximum degree of $G$, the (semi)-saturation number of $F$ with restriction on maximum degree is denoted by
\begin{equation*}
\begin{aligned}
	   (s)sat^{=k}(n,F) & =\min \{ e(G):G  \; \text{is}  \; F\text{-(semi)-saturated, } |G|=n \text{ and } \Delta=k \} \, ; \\
	   (s)sat^{\leq k}(n,F) & =\min \{ e(G):G  \; \text{is}  \; F\text{-(semi)-saturated, } |G|=n \text{ and } \Delta \leq k \} \, ; \\
	   (s)sat^{\geq k}(n,F) & =\min \{ e(G):G  \; \text{is}  \; F\text{-(semi)-saturated, } |G|=n \text{ and } \Delta \geq k \} \, .
\end{aligned}
\end{equation*}
	
	If such an $F$-(semi)-saturated graph with restriction on maximum degree does not exist, we define $(s)sat^{\nabla k}(n,F) =  \infty$ where $\nabla \in \{=,\geq,\leq\}$. For example, $sat^{\nabla r-2} (n, K_r)= \infty$ for $\nabla \in \{= , \leq \}$ and $ n \geq r \geq 3$, since any graph who has maximum degree no more than $r-2$ can not be $K_r$-saturated. Still we have $sat^{\nabla k}(n,F) \geq ssat^{\nabla k}(n,F)$ where $\nabla \in \{=,\geq,\leq\}$. To simplify the notation in our subsequent paper, we would abbreviate $(s)sat^{=k}(n,F)$ as $(s)sat^{k}(n,F)$. For fixed $n,F$ and $ 0 \leq k \leq n-2$, it is clear that
\begin{equation}
(s)sat(n,F)= \min \{ \, (s)sat^{\leq k}(n,F) \, , \; (s)sat^{\geq k+1}(n,F) \, \} \, ,
\end{equation}
and for a special case,
\begin{equation}
(s)sat(n,F)= \min \{ \, (s)sat^{\leq n-2}(n,F) \, , \; (s)sat^{n-1}(n,F) \, \} \, .
\end{equation}
The saturation problem with restriction on maximum degree was first studied by Hajnal \cite{Hajnal} in 1965. F{\"u}redi and Seress \cite{Fur2} found the value of $sat^{\leq \Delta}(n,K_3)$ for $\Delta \geq (n-2)/2$ and sufficiently large $n$ precisely. Later Erd\H os and Holzman \cite{EH} developed the asymptotic value of ${sat^{\leq cn}(n,K_3)}$ as $n \rightarrow \infty$ for $\frac{2}{5} \leq c < \frac{1}{2}$. Alon, Erd\H os, Holzman and Krivelevich \cite{Alon} proved similar results for $K_4$. They also showed that the limit $\lim \limits_{n \rightarrow \infty} \frac{{sat^{\leq cn}(n,K_p)}}{n}$ exists for all $c$ except some sparse value in $(0,1]$.   Amin,  Faudree, Gould and Sidorowicz \cite{Amin} studied $sat^{\leq \Delta}(n,K_p)$ for $\Delta \leq n-2$.
	
	As the  join of graphs plays an elementary role in graph theory, a question is that are we able to determine $sat(n,K_1 \vee F)$ based on the saturation number of $F$\,? Cameron and Puleo \cite{Ca} showed that $sat(n,K_1 \vee F) \leq (n-1)+sat(n-1,F)$ for all $n > |V(F)|$. The reason is that if $G$ is a $(K_1 \vee F)$-saturated graph with a conical vertex $u$, then  $G-u$ is  $F$-saturated. An important problem is to find all graphs such that the equality above holds.
	
\begin{prob}\label{P11} Determine the graph family $\mathcal{F}$ such that for any $F \in \mathcal{F}$, $sat(n,K_1 \vee F) = (n-1)+sat(n-1,F)$ holds for $n$ sufficiently large.
\end{prob}

	Recently, Hu, Luo and Peng \cite{HU} have confirmed that for any graph without isolated vertex but contains an isolated edge, the equality in Problem \ref{P11} holds.
	
\begin{theorem}\label{T12} \text{\bf \cite{HU}} Let $s ,n$ be positive integers and $D$ be a graph without isolated vertex, then for $n \geq 3s^2 -s +2 sat(n-s,D)+1$ we have
\begin{equation*}
\begin{aligned}
sat(n,K_s \vee D) = \binom{s}{2} + s(n-s) + sat(n-s,D) \, .
\end{aligned}
\end{equation*}
\end{theorem}
	
	Although it is not directly stated in Theorem \ref{T12}, the condition that $n \geq 3s^2 -s +2 sat(n-s,D)+1$ holds for large $n$ iff $D$ contains an isolated edge, which is a direct corollary of the results below.

	Let $D$ be a graph and $uv$ be an edge of $D$ such that $d(u) \leq d(v)$. Define the \textit{weight} $wt(uv)$ of the edge $uv$ by $wt(uv)=2 |N(u) \cap N(v)| + |N(v) \setminus N(u)|$. Define the weight of $D$ by $wt(D)=\min \limits_{uv \in E(D)} wt(uv)$. Clearly $wt(D) \geq 1$, and $wt(D) = 1$ iff $D$ contains an isolated edge. Cameron and Puleo \cite{Ca} proved that $sat(n,D) \geq \frac{wt(D)-1}{2}n-\frac{wt(D)^2-4wt(D)+5}{2}$ holds for all $D$ and $n$. Thus if $D$ does not contain an isolated edge, $sat(n,D) \geq \frac{n}{2} - \frac{1}{2}$ holds for all $n \geq wt(D) \geq 2$, which means
\begin{equation*}
3s^2 -s +2 sat(n-s,D)+1 \geq 3s^2 -s +n-s -1+1 > n \,.
\end{equation*}
On the other hand, K\'aszonyi and Tuza \cite{KT} proved that if $D$ contains an isolated edge, then $sat(n,D)=c+o(1)$ for $n$ large, where $c=c(D)$ is a constant only depending on $D$. So for sufficiently large $n$ we have $n \geq 3s^2 -s +2 sat(n-s,D)+1$.

	In order to solve Problem \ref{P11}, we still need to consider the graph $F$ not containing isolated edges. In 2025, Hu, Ji and Cui \cite{HU1} showed that $sat(n,K_1 \vee P_t)=n-1+sat(n-1,P_t)$ holds for $t \geq 5$ and $n$ sufficiently large where $P_t$ is the $t$-vertices path. Song, Hu, Ji and Cui \cite{Song} showed that $sat(n,K_1 \vee C_4)=n-1+sat(n-1,C_4)$ holds for all $n \geq 6$. In this paper we consider the saturation number of wheels, where a wheel graph $W_k=K_1 \vee C_k$. We confirm a slightly stronger result that for all $k \geq 8$, the (semi)-saturation number of wheel graph $W_k$ satisfies the equality of Problem \ref{P11}.
\begin{theorem}\label{T13}
Let $k,n$ be positive integers with $k \geq 8$ and $n \geq 56k^3$, then
\begin{equation*}
    (s)sat(n,W_k)=n-1+(s)sat(n-1,C_k) \, .
\end{equation*}
Moreover, all extremal graphs of $(s)sat(n,W_k)$ must contain a conical vertex while $n,k$ satisfy the condition.
\end{theorem}

	As we consider the (semi)-saturation number of $W_k$, we first need to find out the value of $(s)sat(n,C_k)$. Tuza \cite{TUZ} showed that $(s)sat(n,C_4)= \lfloor \frac{3n-5}{2} \rfloor$ for $n \geq 5$. Chen \cite{Chen1,Chen2} showed that $sat(n, C_5)= \lceil \frac{10(n-1)}{7} \rceil$ for $n \geq 21$. Lan, Shi, Wang and Zhang \cite{LAN} showed that $sat(n,C_6)=\frac{4}{3}n+O(1)$ as $n \geq 9$. For $k \geq 7$, the exact value of $(s)sat(n,C_k)$ is not yet known. The best bounds of $(s)sat(n,C_k)$ are given by F{\"u}redi and Kim \cite{Fur}.

\begin{theorem}\label{T14} \text{\bf \cite{Fur}} (i) For $k \geq 7$ and $n \geq 2k-5$,
\begin{equation*}
\begin{aligned}
    ( 1+\frac{1}{k+2} )n -1    < sat(n,C_k) < (1+\frac{1}{k-4}) n +\binom{k-4}{2} \, .
\end{aligned}
\end{equation*}
(ii) For $n \geq k \geq 6$,
\begin{equation*}
\begin{aligned}
 ( 1+\frac{1}{2k-2} )n -2    < ssat(n,C_k) < (1+\frac{1}{2k-10}) n + k-1 \, .
\end{aligned}
\end{equation*}
\end{theorem}

	In this paper we consider the (semi)-saturation number of $W_k$ with restriction on maximum degree. Our main results are as follow.

\begin{theorem}\label{T15} (i) For $n > k \geq 4$,
\begin{equation*}
\begin{aligned}
    (s)sat^{n-1}(n,W_k)=(n-1) + (s)sat(n-1,C_k) \, .
\end{aligned}
\end{equation*}
(ii) For $ k \geq 6$ and $n \geq 2k^3$,
\begin{equation*}
\begin{aligned}
ssat^{\leq n-2}(n,W_k) \geq \left( \frac{5}{2} - \frac{3}{2k-2} \right)n - \varphi(k)  \, ,
\end{aligned}
\end{equation*}
where $\varphi(k)=(2k^3-17k^2+45k-\frac{57}{2})/(2k-2)$.
\end{theorem}

	By Theorem \ref{T15} and the bounds in Theorem \ref{T14}, we can prove Theorem \ref{T13} immediately.
\vspace{0.5em}

\noindent{\bf Proof of Theorem \ref{T13}. }For $k \geq 8$ and $n \geq 56k^3$, by Theorems \ref{T14} and \ref{T15}, we have
\begin{equation*}
\begin{aligned}
sat^{\leq n-2}(n,W_k) & \geq ssat^{\leq n-2}(n,W_k) \\
& \geq \left( \frac{5}{2} - \frac{3}{2k-2} \right)n - \varphi(k)   \\
& > (n-1)+(1+\frac{1}{k-4}) (n-1) +\binom{k-4}{2} \\
& > (n-1) + sat(n-1,C_k) \\
& =  sat^{n-1}(n,W_k) \geq ssat^{n-1}(n,W_k)\, .
\end{aligned}
\end{equation*}Note that the third inequality holds by $k \geq 8$ and $n \geq 56k^3$. Then by equation (2) we have
\begin{equation*}
\begin{aligned}
 (s)sat(n,W_k) & = \min \{ \, (s)sat^{\leq n-2}(n,W_k) \, , \; (s)sat^{n-1}(n,W_k) \, \} \\
 & = (s)sat^{n-1}(n,W_k)=n-1+(s)sat(n-1,C_k) \, .
\end{aligned}
\end{equation*} 
Clearly, all extremal graphs of $(s)sat(n,W_k)$ must have maximum degree $n-1$ which means all of them would contain a conical vertex.
\vspace{0.5em}

	The remainder of this paper is organized as follow. In next section we would establish several basic results in order to prove the main theorem. Also we would prove Theorem \ref{T15} (i) in this section. In section 3 we would prove Theorem \ref{T15} (ii). We would conclude our paper with some open problems in section 4.

\section{Preliminaries}	
For fixed $k \geq 4$, let $\mathcal{W}$ be the family of all $W_k$-semi-saturated graphs on at least $k+1$ vertices. Let $\mathcal{S}=\{ G \in \mathcal{W} :\text{ there exists } u \in V(G) \text{ such that } d_G(u)=|G|-1 \,\}$ and $\mathcal{T}=\mathcal{W} \setminus \mathcal{S}$. Let $G \in \mathcal{W}$. Since $\delta(W_k)=3$, we have $\delta(G) \geq 2$ and $G$ is connected.
Let $A_G=\{u\in V(G):~d_G(u) \geq k-1\}$, $B_G=\{u\in V(G):~3 \leq d_G(u) \leq k-2\}$ and $C_G=\{u\in V(G):~d_G(u) =2\}$.
We abbreviate a vertex $u \in C_G$ as \textit{root}. The \textit{center} of a given $W_k$ is the vertex of degree $k$. Denote
$$R_G=\{ w \in V(G):\text{$\, \exists \ uv \not \in E(G)$ such that $G+uv$ contains a new $W_k$ centered at $w$} \} .$$
Also we call the vertex in $R_G$ the center of $G$. Clearly $R_G \subseteq A_G$.

\subsection{Properties of roots and distances}
    We first focus on the properties of roots in $G \in \mathcal{W}$ and the distance of two vertices in $G$.
\vspace{0.5em}
	
\begin{lemma}\label{L21}
   For any $u \in V(G)$, $\bigcup \limits_{v \in N[u]} N[v] = V(G)$.
\end{lemma}
\noindent{\bf Proof. }Let $u \in V(G)$.  Since every edges of $W_k$ are in a triangle, we have that a $W_k$-semi-saturated graph must be a $C_3$-semi-saturated graph. So for any $w \in V(G)$, $d(w,u) \leq 2$; otherwise the addition of $wu$ would not create a $C_3$. Thus $\bigcup \limits_{v \in N[u]} N[v] = V(G)$.\q
\vspace{0.5em}

\begin{lemma}\label{L22}
    Let $x,y$ be distinct roots of $G$. Then $xy \not \in E(G)$.
\end{lemma}
\noindent{\bf Proof. }Let $x,y\in C_G$. Suppose for contradiction that $N(x)=\{y,w\}$ for some $w$. Pick a vertex $z$ distinct from $y,w$. Since the addition of $xz$ would create a $W_k$ in $G+xz$, such a $W_k$ must contain $x,y,z,w$ which contradicts $d(y)=2$. \q
\vspace{0.3em}

\begin{lemma}\label{L23}
    Let $x,y$ be distinct roots of $G$. Then $|N(x) \cap N(y)|=1$.
\end{lemma}
\noindent{\bf Proof. }Let $x,y\in C_G$. Assume that $N(x)=\{w,z\}$. By Lemma \ref{L22}, $xy\notin E(G)$. Then the addition of $xy$ would create a $W_k$ in $G+xy$ containing $x,y,w,z$. Since $R_G \subseteq A_G$ and $x,y \in C_G$ we can assume that $w \in R_G$. Then $w \in N(y)$ and $z \not \in N(y)$; otherwise such a $W_k$ is $W_3$, a contradiction to $k \geq 4$. \q
\vspace{0.3em}

	If $C_G$ is not empty, let $C_G=\{x_1,\ldots,x_r\}$ for some positive $r \in \mathbb{N}$. By Lemma \ref{L23}, either there exists some $u$ such that $N(x_i)=\{u,y_i\}$ for all $i \in [r]$ where $y_i \neq y_j$ for $1 \leq i < j \leq r$, or $r=3$ and $N(x_1)=\{y,z\}$, $N(x_2)=\{z,w\}$, $N(x_3)=\{w,y\}$ for some distinct $y,z,w$. Anyway we know that $|N(C_G)| \geq |C_G|$.

\subsection{Properties on $\mathcal{S}$}

For a graph $G$, recall that $u \in V(G)$ is a conical vertex if $d_G(u)=|G|-1$. By our definition, every $G \in \mathcal{S}$ has a conical vertex.

\begin{prop}\label{P28}
For any $k \geq 4$ and $n \geq k+1$, let $G$ be a $W_k$-(semi)-saturated graph with $|G|=n$. Assume that $u $ is a conical vertex of $G$, then $G-u$ is $C_k$-(semi)-saturated.
\end{prop}
\noindent{\bf Proof. }Let $H=G-u$. If $G$ is $W_k$-free, then $H$ is $C_k$-free since $d(u)=n-1$. Let $xy$ be a nonedge in $H$. Then the addition of $xy$ in $G$ would create a new $W_k$ (says $W^*$) centered at some $w$.

    If $u \not \in V(W^*)$, then $W^*$ is a subgraph of $H+xy$ and the addition of $xy$ would create new $C_k$ in $H$. If $u \in V(W^*)$ and $u=w$, then $W^* -u$ is a $C_k$ in $H+xy$ containing $xy$. If $u \in V(W^*)$ but $u \neq w$, then $u,w$ are both vertices of maximum degree $k$ in $W^*$. We can exchange $w$ and $u$ in $W^*$ and are able to find a $C_k$ in $H+xy$. Thus $H$ is $C_k$-(semi)-saturated. \q
\vspace{0.5em}

\noindent {\bf Proof of Theorem \ref{T15} (i). }Let $G$ be a $W_k$-(semi)-saturated graph of order $n$ with minimum number of edges and $u$ be the conical vertex of $G$. By Proposition \ref{P28}, $H=G-u$ is $C_k$-(semi)-saturated, then
\begin{equation}
(s)sat^{n-1}(n,W_k) = e(G) = n-1 + e(H) \geq (n-1) + (s)sat(n-1,C_k) .
\end{equation}
On the other hand, let $H$ be a $C_k$-(semi)-saturated graph such that $e(H)= (s)sat(n-1,C_k)$. Then $H\vee K_1$ is $W_k$-(semi)-saturated. Thus the equality in (3) holds and the proof is done. \q

\section{The family $\mathcal{T}$}

Let $G \in \mathcal{T}$. In this section we establish the lower bound of size of $G$. Fixed $k \geq 6$ and $n \geq (k-1)^2+1 $, assume that $G \in \mathcal{T}$ with $|G|=n$. Then $\Delta(G)\le n-2$. Let $A^*_G=B_G \cup C_G$, we begin from an easy lemma.

\begin{lemma}\label{L31}
    Let $r \in \mathbb{N}$. If $|A^*_G| \leq r$, then $e(G) \geq \frac{(k-1)n}{2}-\frac{(k-3)r}{2}$.
\end{lemma}
\noindent{\bf Proof. } We have
\begin{equation*}
\begin{aligned}
e(G) & = \frac{1}{2} \sum_{v \in V(G)} d(v) = \frac{1}{2} \left( \sum_{v \in A_G} d(v) + \sum_{v \in A^*_G} d(v)  \right) \\
& \geq \frac{1}{2} \left( (k-1)(n-r) + 2r  \right) = \frac{(k-1)n}{2}-\frac{(k-3)r}{2} \, .
\end{aligned}
\end{equation*}\q
\vspace{0.5em}

By Lemma \ref{L31},	we can assume that $|A_G^*| \geq k$, else $e(G) \geq \frac{5n}{2}-\frac{k(k-3)}{2}$ which immediately leads to our conclusion. For any $u \in A^*_G$, we have $ N(u) \cap A_G \neq \emptyset$ since there must exists a $v  \in A_G^* \setminus N(u)$ and the addition of $uv$ would create a $W_k$ centered at some $w \in R_G \subseteq A_G$.

\begin{lemma}\label{L32}
    If $\bigcap \limits_{v \in A_G^*} N(v)= \emptyset$, then there exists $X \subseteq V(G)$ with $|X| \leq (k-1)^2$, such that $e(u,X) \geq 2$ for any $ u \in V(G)\setminus X$.
\end{lemma}
\noindent{\bf Proof. }Let $x \in A_G^*$ and $N(x) \cap A_G=\{y_1, \ldots, y_\ell \}$. Then $\ell \leq k-2$. Since $\bigcap \limits_{v \in A_G^*} N(v)= \emptyset$, for any $y_i$ ($1\le i\le \ell$), there exists $z_i \in A_G^*$ such that $y_iz_i \not \in E(G)$. Let $X= \left( \bigcup \limits_{i=1}^\ell N[z_i] \right) \cup N[x]$. Then $|X| \leq (k-1)^2$. Since $n \geq (k-1)^2+1$, there exists $u \in V(G)\setminus X$. Then $ux,uz_i \not \in E(G)$ for any $1 \leq i \leq \ell$.

	The addition of $ux$ in $G$ would create a $W_k$. If the center of such a $W_k$ is $u$, then $|N(u) \cap N(x)| \, \geq 2$, which means $e(u,X) \geq 2$. Assume the center of such a $W_k$ is not $u$ but some $y_i$ (says $y_1$). Since the addition of $uz_1$ would create a $W_k$ not centered at $y_1$ as $y_1z_1 \not \in E(G)$,  there exists $w \in N(z_1) \cap N(u)$ distinct from $y_1$. Then $w,y_1 \in N(u) \cap X$, and hence $e(u,X) \geq 2$. \q
\vspace{0.3em}

\begin{lemma}\label{L33}
    Let $x \in C_G$ and $N(x)=\{y,z\}$, then $d(y),d(z) \geq k$. Moreover, $yz \in E(G)$.
\end{lemma}
\noindent{\bf Proof. }By symmetry we only need to prove $d(z) \geq k$. Since $G \in \mathcal{T}$ there exists $w \not \in N(y)$. The addition of $xw$ can only create a $W_k$ containing $x,y,z,w$ which is centered at $z$ since $yw \not \in E(G)$. So $d(z) \geq k$ and $yz \in E(G)$. \q
\vspace{0.5em}

\begin{lemma}\label{L34}
    For a vertex set $R \subseteq V(G)$, we have $\sum \limits_{v \in R} d(v) \geq 3|R| + \varepsilon (R)$, where $$ \varepsilon(R) =  \max \{ -|C_G|, -(n-|R|) \} . $$
\end{lemma}
\noindent{\bf Proof. }It is clear that $\sum \limits_{v \in R} d(v) \geq 3|R|-  |C_G \cap R|$. By Lemma \ref{L33} and $|N(C_G)| \geq |C_G|$, there are at least $ \max \{0 , |C_G|-(n-|R|) \}$ vertices in $R$ with degree at least $k$. So
\begin{equation*}
\begin{aligned}
\sum \limits_{v \in R} d(v) & \geq 3|R|-  |C_G | + (k-3) \max \{ 0 , |C_G|-(n-|R|)\} \\
 & \geq  3|R|+ \max \{ -|C_G|, -(n-|R|) \} \, .
\end{aligned}
\end{equation*}
\q
\vspace{0.5em}

\begin{prop}\label{P35}
 If $\bigcap \limits_{v \in A_G^*} N(v)= \emptyset$, then $e(G) \geq \frac{5}{2}n- 3(k-1)^2$.
\end{prop}
\noindent{\bf Proof. }We take the set $X$ as Lemma \ref{L32} stated. Then $\sum \limits_{v \in X} d(v) \geq e(X,V(G)\setminus X) \geq 2(n-|X|)$. By Lemma \ref{L34}, $\sum \limits_{v \in V(G)\setminus X} d(v) \geq 3(n-|X|)+\varepsilon(V(G)\setminus X)$, where
 $$ \varepsilon(V(G)\setminus X) =   \max \{ -|C_G|, -|X| \}   \geq -(k-1)^2. $$
So we have
\begin{equation*}
\begin{aligned}
e(G) & = \frac{1}{2} \sum \limits_{v \in V(G)} d(v) = \frac{1}{2} \left( \sum \limits_{v \in X} d(v) + \sum \limits_{v \in V(G)- X} d(v) \right) \\
& \geq \frac{1}{2} \left(  2(n-|X|) + 3(n-|X|) - (k-1)^2 \right) \\
& \geq \frac{5}{2}n- 3(k-1)^2,
\end{aligned}
\end{equation*}and we are done. \q
\vspace{0.5em}

	Proposition \ref{P35} establish the lower bound of $e(G)$ when $\bigcap \limits_{v \in A_G^*} N(v)= \emptyset$. In the next subsection we would consider the case $\bigcap \limits_{v \in A_G^*} N(v)\not= \emptyset$. Let $a^* \in \bigcap \limits_{v \in A_G^*} N(v) \subseteq A_G$. We call the common neighbour $a^*$ the \textit{grand connector}. For $x \in C_G$, let $N(x)=\{a^*,y\}$. We say that $y$ is the \textit{related connector} of $x$ and $x$ is the \textit{related root} of $y$. The set of all related connectors is denoted by $C^*$. Clearly $|C^*|=|C_G|$ by Lemma \ref{L23}.

\subsection{Graph with a grand connector}

We fixed the grand connector $a^* \in \bigcap \limits_{v \in A_G^*} N(v) \subseteq A_G$. Since $G \in \mathcal{T}$, there exists $b^* \in A_G$ such that $a^*b^* \not \in E(G)$. By Lemma \ref{L33}, $A_G^* \cup C^* \subseteq N(a^*)$. Thus $d(a^*) \geq |A_G^*|+ |C^*| = |B_G|+2|C_G|$.
\vspace{0.5em}

\begin{definition}\label{D36}
	We define three functions $f_i$ on $G$. Let
\begin{equation*}
\begin{aligned}
f_1(G) & =  2|A_G|+\frac{5}{2}|B_G|+\frac{5}{2}|C_G|- d(b^*)-\frac{7}{2} \, ; \\
f_2(G) & = \frac{(k-1)}{2} |A_G| + 2|B_G| + \frac{5}{2} | C_G|+\frac{1}{2} d(b^*) -k +1  \, ; \\
f_3(G) & = |A_G| + \frac{5}{2} |B_G| + 3 |C_G| - (k^2-3k+3) \, .
\end{aligned}
\end{equation*}
\end{definition}

	We would prove that $e(G) \geq \max \{ f_1(G), f_2(G), f_3(G) \}$.
\vspace{0.5em}

\begin{lemma}\label{L37}
    $e(G) \geq f_1(G)$.
\end{lemma}
\noindent{\bf Proof. }Let $ R=V(G)- (N[b^*] \cup \{a^*\})$, apply Lemma \ref{L34} to $R$ we have $\sum \limits_{v \in R} d(v) \geq 3|R| + \varepsilon (R)$, where $\varepsilon(R)=   \max \{ -|C_G|, -(d(b^*)+2) \}  \geq -|C_G|$. By Lemma \ref{L21} $\bigcup \limits_{v \in N[b^*]} N[v] = V(G)$, hence $e(N(b^*), V(G) - N[b^*]) \geq n-d(b^*)-1$ and 
\begin{equation*}
\sum \limits_{v \in N(b^*)} d(v) \geq e(N(b^*), V(G) - N[b^*]) + e(N(b^*),b^*) \geq n-1 \, .
\end{equation*}
So we have 
\begin{equation*}
\begin{aligned}
2e(G) & = \sum \limits_{v \in V(G)} d(v) = d(a^*)+d(b^*)+ \sum \limits_{v \in R} d(v) + \sum \limits_{v \in N(b^*)} d(v) \\
& \geq |B_G|+2|C_G|+d(b^*) + 3(n-d(b^*)-2) + \varepsilon (R) + n-1 \\
& \geq 4|A_G|+5|B_G|+5|C_G|-2d(b^*)-7 = 2f_1(G),
\end{aligned}
\end{equation*}
and we are done. \q
\vspace{0.5em}

\begin{lemma}\label{L38}
	$e(G) \geq f_2(G)$.
\end{lemma}
\noindent{\bf Proof. }By Lemma \ref{L33}, every vertex in $C^* \subseteq A_G$ has degree at least $k$. Since $a^*b^* \not \in E(G)$, we have $b^* \in A_G \setminus C^*$. Then
\begin{equation*}
\begin{aligned}
2e(G) & = \sum \limits_{v \in V(G)} d(v) = d(a^*)+d(b^*)+ \sum \limits_{v \in B_G} d(v) + \sum \limits_{v \in C_G} d(v) + \sum \limits_{v \in A_G- \{a^*,b^*\}} d(v) \\
& \geq |B_G|+2|C_G|+d(b^*)+3|B_G| + 2|C_G| + (k-1)(|A_G|-2)+|C^*| \\
& =  (k-1) |A_G| + 4|B_G|+5|C_G| + d(b^*) -2k+2 = 2f_2(G) \, ,
\end{aligned}
\end{equation*}
and we are done. \q
\vspace{0.5em}

\begin{lemma}\label{L39}
	$e(G) \geq f_3(G)$.
\end{lemma}
\noindent{\bf Proof. }We define the middle-prioritized weight $wt(u)$ on every vertex $u$ according to edges as follows. For an edge $uv\in E(G)$\,: \\
\indent (1) assign $1$ to $u$ and $0$ to $v$ if $u \in A_G^*$ and $v \in A_G$\,; \\
\indent (2) assign $\frac{1}{2}$ to $u$ and $\frac{1}{2}$ to $v$ if $u,v \in B_G$\,; \\
\indent (3) assign $0$ to $u$ and $0$ to $v$ if $u,v \in A_G$\,. \\
\noindent Then $wt(u)$ is the total obtained value. By Lemma \ref{L22}, $C_G$ is an independent set. By Lemma \ref{L33} and the definition of $a^*$, $e(B_G,C_G)=0$. So every edge in $E(G) - E(G[A_G])$ has assigned total weight 1 to its endpoints, and every edge in $ E(G[A_G])$ has assigned total weight 0 to its endpoints. Then $ e(G) =  \sum \limits_{v \in V(G)} wt(v) + e(G[A_G])$. By Lemma \ref{L33}, $wt(u)=2$ for $u \in C_G$. Also $wt(v) =0$ for $v \in A_G$ and $wt(w) \geq 2$ for $w \in B_G$ as $a^* \in N(w)\cap A_G$.

	Partition $B_G$ into $B_G^1 = \{ v \in B_G: wt(v)=2 \}$ and $B_G^2=B_G \setminus B_G^1$. Then for any $v \in B_G^1$, $N(v) \cap A_G=\{a^*\}$ and $|N(v) \cap B_G|=2$ which implies $d(v)=3$ (combining with $e(B_G,C_G)=0$).
\vspace{0.3em}

\noindent {\bf Claim 1. }For $x \in B_G$, if $|N(x) \cap B_G^1| \geq 1$, then $xb^* \in E(G)$. Moreover, we have $x \in B_G^2$. As a corollary, $B_G^1$ is an independent set and $|N(y) \cap B_G^2|=2$ for all $y \in B_G^1$.

\noindent {\bf Proof of Claim 1. }Pick $y \in N(x) \cap B_G^1$ and suppose that $N(y)=\{a^*,x,z\}$. Then $z\in B_G$. The addition of $yb^*$ would create a $W_k$ which is centered at $b^*$ and contains $x,z$ by $a^*b^* \not \in E(G)$. Thus $xb^* \in E(G)$. Since $a^*x\in E(G)$, we have $wt(x)\ge 5/2$ and then $x \not \in B_G^1$. \qed
\vspace{0.5em}

\noindent {\bf Claim 2. }For $x \in B_G$, if $|N(x) \cap B_G^1| =t \geq 1$, then $wt(x) \geq \max \{3, 2+ \frac{1}{2}t\}$.

\noindent {\bf Proof of Claim 2. }By Claim 1  and the definition of $a^*$, we have $a^*,b^* \in N(x)$. If $t \geq 2$ then we are done, so we only need to show that $wt(x) \geq 3$ for $t=1$. Suppose for contradiction that $wt(x)<3$ for $t=1$. Then $d(x)=3$. Let
$N(x)=\{a^*,b^*,y\}$ for some $y \in B_G^1$. We can assume $N(y)=\{a^*,x,z\}$ where $z \in B_G^2$. Then $xz \not \in E(G)$. The addition of $xz$ would create a $W_k$ which is centered at $a^*$ and does not contain $b^*$ by $N(x) \cap N(b^*)=\emptyset$. Moreover, such a $W_k$ must contain $a^*,x,y,z$ which means $k=3$, a contradiction. \qed
\vspace{0.5em}

\noindent {\bf Claim 3. }$\sum \limits_{v \in B_G} wt(v) \geq \frac{5}{2} |B_G|$.

\noindent {\bf Proof of Claim 3. }Let $\mathcal{H}_1=(V_1,E_1)$ be a subgraph of $G$, where $V_1=B_G$ and $E_1=\{uv \in E(G): u \in B_G^1, v \in B_G^2\}$.  Let $H_1,\ldots, H_r$ be the components of $\mathcal{H}_1$. We prove a slightly stronger results that $\sum \limits_{v \in V(H_i)} wt(v) \geq \frac{5}{2} |V(H_i)|$ holds for all $i \in [r]$.

	By symmetry we only need to consider $H_1$. If $H_1$ is an isolated vertex $x$, by Claim 1, $ x \in B_G^2$. Hence $\sum \limits_{v \in V(H_i)} wt(v) = wt(x) \geq \frac{5}{2}$. Now we assume that $|V(H_1)|\ge 2$. Then $|V(H_1)\cap B_G^1|\ge 1$ and $|V(H_1)\cap B_G^2|\ge 2$ by Claim 1. Let $V(H_i)=X\cup Y$, where   $X=\{x_1 ,\ldots ,x_q \}\subseteq B_G^2$ with $q \geq 2$ and $Y=\{y_1, \ldots, y_p\}\subseteq B_G^1$ with $p \geq 1$. Define $r_i=|\{ x \in X: |N(x) \cap Y|=i \} |$. Then $\sum \limits_{i=1}^p r_i = q$. Since $H_1$ is a bipartite graph, we have $2p = \sum \limits_{i=1}^p ir_i$. Then by Claim 2,
\begin{equation*}
\begin{aligned}
\sum \limits_{v \in H_1} wt(v) & = \sum \limits_{v \in Y} wt(v) + \sum \limits_{v \in X} wt(v) \\
& = 2p + \sum \limits_{v \in X: \, |N(x) \cap Y|=1} wt(v) + \sum \limits_{v \in X: \,|N(x) \cap Y| \geq 2} wt(v) \\
& \geq 2p +3r_1 +  \sum \limits_{i=2}^p (2+\frac{1}{2}i) r_i
 = \frac{5}{2}p + \frac{11}{4}r_1+ \sum \limits_{i=2}^p (2+\frac{1}{4}i) r_i \\
& \geq \frac{5}{2}p + \frac{5}{2} \sum \limits_{i=1}^p r_i = \frac{5}{2} (p+q) = \frac{5}{2} |V(H_1)|.
\end{aligned}
\end{equation*}
So the proof is done. \qed
\vspace{0.5em}

	In order to compute $e(G)$, we still need to bound $ e(G[A_G])$.
\vspace{0.5em}

\noindent {\bf Claim 4. }The number of  components of $G[A_G]$ is at most $(k-2)(k-1)$.

\noindent {\bf Proof of Claim 4. }Let $G[A_G] = \mathcal{H}_2=(V_2,E_2)$ where $V_2 = A_G$. Suppose that $\mathcal{H}_2$ contains $r$  components $H_1,\ldots, H_r$. If $A_G^* = \emptyset$ then $r=1$ by $G $ being connected. Assume $A_G^* \not= \emptyset$. Let $x \in A_G^*$. By Lemma \ref{L21}, $\bigcup \limits_{v \in N[x]} N[v] = V(G)$. Then for all $ i \in [r]$ one of the following must hold\,:

\indent (1) $\exists \, y \in V(H_i)$ such that $xy \in E(G)$; \\
\indent (2) $\exists \, y \in V(H_i)$ and $z \in N(x) \cap A_G^*$ such that $yz \in E(G)$.

\noindent Since $d(x) \leq k-2$, at most $k-2$ components of  $\mathcal{H}_2$ satisfy (1). Since any vertex in $N(x) \cap A_G^*$ has degree no more than $k-2$, at most $(k-2)^2$ components of  $\mathcal{H}_2$ satisfy (2). Thus we have $r \leq (k-2)^2 + (k-2) = (k-2)(k-1)$. \qed
\vspace{0.5em}

\noindent {\bf Claim 5. }$e(G[A_G]) \geq |A_G|+|C_G|-(k^2-3k+3)$.

\noindent {\bf Proof of Claim 5. }Since $a^*b^* \not \in E(G)$, $b^* \not \in C^* \cup \{a^*\}$. By Lemma \ref{L21},  $\bigcup \limits_{v \in N[b^*]} N[v] = V(G)$. Since $ C_G\subseteq \bigcup \limits_{v \in N[b^*]} N[v]$, we have $xb^*\in E(G)$ for any $x\in C^*$. Since $\{a^*,b^*\} \cup C^* \subseteq A_G$, there are at least $|C^*|-1 = |C_G|-1$ cycles in $G[A_G]$. Thus
\begin{equation*}
\begin{aligned}
e(G[A_G]) & \geq |A_G|- r + c \geq |A_G|+|C_G|-(k^2-3k+3) \, ,
\end{aligned}
\end{equation*}
where $r \leq (k-2)(k-1)$ is the number of  components in $G[A_G]$ and $c \geq |C_G|-1$ is the number of
cycles in $G[A_G]$. \qed
\vspace{0.5em}

	By Claims 3 and 5 we have
\begin{equation*}
\begin{aligned}
e(G)&  =  \sum \limits_{v \in V(G)} wt(v) + e(G[A_G]) \\
& = \sum \limits_{v \in C_G} wt(v) + \sum \limits_{v \in B_G} wt(v)  + e(G[A_G]) \\
& \geq 2|C_G| + \frac{5}{2} |B_G|  + |A_G|+|C_G|-(k^2-3k+3) \\
& = 3 |C_G| + \frac{5}{2} |B_G| + |A_G| - (k^2-3k+3) = f_3(G) \, .
\end{aligned}
\end{equation*}
The proof is done. \q
\vspace{0.5em}

\begin{prop}\label{P310}
 For $k \geq 6$, if $\bigcap \limits_{v \in A_G^*} N(v) \neq \emptyset$,
then $e(G) \geq \left( \frac{5}{2} - \frac{3}{2k-2} \right)n - \varphi(k)$ where $\varphi(k)=(2k^3-17k^2+45k-\frac{57}{2})/(2k-2)$.
\end{prop}
\noindent{\bf Proof. }By the discussion above in this section and Lemmas \ref{L37}, \ref{L38}, \ref{L39}, we have $e(G) \geq \max \{f_1(G), f_2(G), f_3(G) \}$. Hence (recall $n=|A_G|+|B_G|+|C_G|$)\,: 
\begin{equation*}
\begin{aligned}
e(G)&  \geq \frac{3f_1(G)+6f_2(G)+(2k-11)f_3(G)}{2k-2} \\
& = (\frac{5k-8}{2k-2})(|A_G|+|B_G|+|C_G|) + (\frac{2k-5}{4k-4}) |C_G| - (\frac{2k^3-17k^2+45k-\frac{57}{2}}{2k-2}) \\
& \geq \left( \frac{5}{2} - \frac{3}{2k-2} \right)n - \varphi(k) \, .
\end{aligned}
\end{equation*}
The proof is done. \q

\subsection{On the lower bound of $ssat^{\leq n-2}(n,W_k)$}

We are now able to prove Theorem \ref{T15} (ii).
\vspace{0.5em}
	
\noindent {\bf Proof of Theorem \ref{T15} (ii). }For $k \geq 6$ and $n \geq 2k^3$, let $G \in \mathcal{S}$ of order $n$ with minimum number of edges. If $\bigcap \limits_{v \in A_G^*} N(v)= \emptyset$, then by Proposition \ref{P35} $e(G) \geq \frac{5}{2}n- 3(k-1)^2=:h_1$. Otherwise by Proposition \ref{P310} $e(G) \geq \left( \frac{5}{2} - \frac{3}{2k-2} \right)n - \varphi(k)=:h_2$ where $\varphi(k)=(2k^3-17k^2+45k-\frac{57}{2})/(2k-2)$. Hence
\begin{equation*}
ssat^{\leq n-2}(n,W_k) = e(G) \geq \min\{ h_1,h_2 \} = h_2 = \left( \frac{5}{2} - \frac{3}{2k-2} \right)n - \varphi(k) .
\end{equation*} \q

\section{Concluding remark}

    We have proved all Theorems in this paper by the discussion above. Although we only show that $W_k \in \mathcal{F}$ where $k \geq 8$ and $\mathcal{F}$ is the graph family defined in Problem \ref{P11}, we believe that such result is correct for all $ k \geq 4$.

\begin{conj}\label{C41}
For $k \geq 4$ there exists a constant $N_k$ such that for all $n>N_k$ we have $(s)sat(n,W_k)=n-1+(s)sat(n-1,C_k)$. Moreover, if $k \geq 6$ then the extremal graph of $sat(n,W_k)$ must contain a conical vertex. 
\end{conj}

	Song, Hu, Ji and Cui \cite{Song} confirmed that $W_4 \in \mathcal{F}$ but the extremal graph of $sat(n,W_4)$ may not contain a conical vertex. Due to the complex constructions of $sat(n,C_5)$ we speculate that there exists some extremal graph of $sat(n,W_5)$ not containing conical vertex. For $k \geq 6$ we believe that all extremal graphs of $sat(n,W_k)$ must contain a conical vertex.

    In order to determine the exact value of $(s)sat(n,W_k)$ for $ k \geq 8$ and $n $ sufficiently large, our Theorem shows that it is inevitable to find the exact value of $(s)sat(n,C_k)$. This is a vital and intriguing problem.

\begin{prob}\label{P42}
Determine the exact value of $(s)sat(n,C_k)$ for $ k \geq 7$ and $n$ sufficiently large.
\end{prob}

\section*{Acknowledgement}
The research is supported by Beijing Natural Science Foundation (No.~1244047); and by the National Natural Science Foundation of China (No.~12171272, 12401445), China Postdoctoral Science Foundation (No.~2023M740207).

\end{document}